\documentclass[12pt,a4paper,reqno]{article}%
\usepackage{amssymb}
\usepackage{amsmath}
\usepackage{amsfonts}
\usepackage{graphicx}
\usepackage{rotating}
\usepackage{mathrsfs}
\usepackage[all]{xy}
\usepackage{hyperref}%
\providecommand{\U}[1]{\protect\rule{.1in}{.1in}}
\SelectTips{cm}{}
\newtheorem{theorem}{Theorem}

\newtheorem{problem}{Problem}

\newtheorem{proposition}[theorem]{Proposition}

\begin{document}

\title{\textbf{Geometric Hamilton-Jacobi Field Theory}}
\author{\textsc{{L.{} Vitagliano}\thanks{\textbf{e}-\textit{mail}:
\texttt{lvitagliano@unisa.it}}} \\ {} \\ {\small {DipMat, Universit\`a degli Studi di Salerno, and}} \\ {\small {Istituto Nazionale di Fisica Nucleare, GC Salerno}} \\ {\small {Via Ponte don Melillo, 84084 Fisciano (SA), Italy}} \\ {} \\ {\small {Istituto Tullio Levi-Civita}} \\ {\small {via Colacurcio 54, 83050 Santo Stefano del Sole (AV), Italy}} }
\maketitle


\begin{abstract}
I briefly review my proposal about how to extend the geometric Hamilton-Jacobi theory to higher derivative field theories on fiber bundles. 
\end{abstract}

\newpage

\section*{Introduction}

The standard Hamilton-Jacobi (HJ) problem consists in finding Lagrangian submanifolds of the phase space which are preserved by the dynamics. In its turn, the HJ problem is locally equivalent to the HJ Equation. The restricted dynamics on a solution of the HJ problem is simpler than the complete one and, therefore, solving the HJ problem helps in finding solutions of the Hamilton equations. In the series of papers \cite{c...06,c...09,c...10}, Cari\~{n}ena, Gr\`{a}cia, Marmo, Mart\'{\i}nez, Mu\~{n}oz-Lecanda, Rom\'{a}n-Roy, presented a generalized HJ problem on both the Lagrangian and the Hamiltonian side of a Lagrangian theory. Their formulation \textquotedblleft is based on the idea of obtaining solutions of a second order differential equation by lifting solutions of an adequate first order differential equation\textquotedblright, and I will refer to it as \emph{Geometric HJ Theory}.

Inspired by the above mentioned papers \cite{c...06,c...09,c...10}, I recently proposed an extention of the geometric HJ theory to higher derivative field theories on fiber bundles \cite{v10,v11} (see also \cite{dL...08} for the first order case). My proposal is based on the Lagrangian-Hamiltonian formalism developed in \cite{v10b} (see also \cite{c...09b}) and is manifestly coordinate independent and natural, i.e., independent of geometric structures other than the action functional. I will refer to my extension of the geometric HJ theory as \emph{Geometric HJ Field Theory}.

The aim of the paper is to review the geometric HJ field theory. The main results will be presented without proofs. We refer to \cite{v10,v11} for details.

\section*{Geometric Hamilton-Jacobi Theory}

I first review the geometric HJ theory. For simplicity I only consider the autonomous case. I present an equivalent formulation on the fibered product of the velocity space and the phase space. Such formulation is suitable for the extension to field theory.

Let $Q$ be a configuration manifold, and $TQ$ (resp.{} $T^\ast Q$) the corresponding velocity (resp.{} phase) space. A \emph{Lagrangian theory} on $Q$ is the datum of a Lagrangian function $L$ on $TQ$. It determines a system of (second order) Euler-Lagrange (EL) equations $\mathscr{E}_{EL} \subset T^2 Q$ for curves $\gamma : \mathbb{I} \longrightarrow Q$ in $Q$. Roughly speaking, extremals of the action functional $\int L dt$ are solutions of the EL equations, i.e., curves $\gamma$ in $Q$ such that their second lift $\ddot{\gamma} : \mathbb{I} \longrightarrow T^2 Q$ take values in $\mathscr{E}_{EL}$. Furthermore, $L$ determines a system of first oder, implicit Hamilton equations $\mathscr{E}_H \subset TT^\ast Q$ for curves in $T^\ast Q$. There are several equivalent definitions of $\mathscr{E}_H$. A possible one is the following.    

Put $\mathbb{T}Q := TQ \times_Q T^\ast Q$, and let $p : \mathbb{T}Q \longrightarrow T^\ast Q$ and $q : \mathbb{T}Q \longrightarrow TQ$ be canonical projections. First of all the Lagrangian $L$ determines a first order ODE for curves in $\mathbb{T}Q$ as follows. Namely, let $E_0 : \mathbb{T}Q \longrightarrow \mathbb{R}$ be the function determined by the canonical pairing between $TQ$ and $T^\ast Q$, i.e., $E_0 (\xi, \theta) := \theta (\xi)$, $(\theta, \xi) \in \mathbb{T}Q$. Put $E := E_0 - q^\ast (L)$. Moreover, let $\Omega$ be the pull-back of the canonical symplectic structure on $T^\ast Q$ via $p$. We can consider the following equation for a curve $\Sigma$ in $\mathbb{T}Q$:
\begin{equation}
i({\dot{\Sigma}}) \Omega = dE |_\Sigma \label{ELHODE}.
\end{equation}
As usual, the geometric portrait of Eq.{} (\ref{ELHODE}) is a subset $\mathscr{E}_{ELH}$ of $T\mathbb{T}Q$. Since $\mathbb{T}Q$ projects onto $T^\ast Q$, then $\mathscr{E}_{ELH}$ projects onto a suitable subset $\mathscr{E}_H$ of $TT^\ast Q$ which, in its turn, can be interpreted as an implicit first order ODE for curves in $T^\ast Q$. $\mathscr{E}_H$ is precisely the implicit Hamilton Equations.

Eq.{} (\ref{ELHODE}) covers both the EL Lagrange and the implicit Hamilton Equations, i.e., if $\Sigma$ is a solution of Eq.{} (\ref{ELHODE}), then $q \circ \Sigma = \dot{\gamma}$ for a solution $\gamma$ of the EL Equations, and $p \circ \Sigma$ is a solution of the implicit Hamilton Equations. For this reason, I will refer to Eq.{} (\ref{ELHODE}) as the \emph{Euler-Lagrange-Hamilton (ELH) ODE}. Solutions of the ELH ODE take values in the graph $\mathscr{P} \subset \mathbb{T}Q$ of the Legendre transform $FL: TQ \longrightarrow T^\ast Q$. Notice that the canonical projection $\mathscr{P} \longrightarrow Q$ is a bundle whose sections are pairs $(\nabla, T)$ such that $\nabla$ is a vector field on $Q$, and $T = FL \circ X$.

Let $\ldots, q^i, \dots$ be coordinates on $Q$, and $\ldots, \dot{q}{}^i, \ldots$ (resp.{}, $\ldots, p_i, \ldots$) associated fiber coordinates on $TQ$ (resp., $T^\ast Q$). Recall that the EL Equations have the local expression:
\begin{equation*}
\frac{\partial L}{\partial q^i}\left(q,\tfrac{dq}{dt}\right) - \frac{d}{dt} \frac{\partial L}{\partial \dot{q}{}^i}\left(q,\tfrac{dq}{dt}\right) =0.
\end{equation*}
The ELH ODE has the following coordinate expression
\begin{equation*}
\left\{ 
\begin{array}{l}
\dfrac{d}{dt}q^{i}=\dot{q}{}^{i} \\ 
\vspace{-0.2cm} \\ 
\dfrac{d}{dt}p_{i}=\dfrac{\partial L}{\partial q^{i}}\left( q,\dot{q}\right) 
\\ 
\vspace{-0.2cm} \\ 
p_{i}=\dfrac{\partial L}{\partial \dot{q}{}^{i}}\left( q,\dot{q}\right) 
\end{array}.
\right.
\end{equation*} 
In particular $\mathscr{P}$ is locally given by
\begin{equation*}
p_{i}=\dfrac{\partial L}{\partial \dot{q}{}^{i}}\left( q,\dot{q}\right).
\end{equation*}
The implicit Hamilton Equations have the local expression
\begin{equation*}
\left\{ 
\begin{array}{l} 
\dfrac{d}{dt}p_{i}=\dfrac{\partial L}{\partial q^{i}}\left( q,\tfrac{dq}{dt}\right) 
\\ 
\vspace{-0.2cm} \\ 
p_{i}=\dfrac{\partial L}{\partial \dot{q}{}^{i}}\left( q,\tfrac{dq}{dt}\right) 
\end{array}.
\right.
\end{equation*}

The central problem in geometric HJ theory has got the following equivalent formulations

\begin{problem}[Lagrangian Generalized HJ Problem]
Find a vector field $\nabla$ on $Q$ such that any integral curve of $\nabla$ is a solution of the EL Equations
\end{problem}

\begin{problem}[Hamiltonian Generalized HJ Problem]
Find a vector field $\nabla$ on $Q$ such that, for any integral curve $\gamma$ of $\nabla$, $FL \circ \dot{\gamma}$ is a solution of the implicit Hamilton Equations.
\end{problem}

The key result in the geometric HJ theory is the following

\begin{theorem}[Generalized HJ Theorem]
Let $\nabla$ be a vector field on $Q$. The following conditions are equivalent:
\begin{itemize}
\item $\nabla$ is a solution of the generalized HJ problem;
\item $i(\nabla) (\nabla , T)^\ast \Omega = (\nabla, T)^\ast E$, with $(\nabla , T)$ a section of $\mathscr{P} \longrightarrow Q$.
\end{itemize}
\end{theorem}
Therefore, solving the generalized HJ problem amounts to solving the PDE 
\begin{equation}
i(\nabla) (\nabla , T)^\ast \Omega = (\nabla, T)^\ast E, \label{GHJE}
\end{equation}
for sections $(\nabla, T)$ of $\mathscr{P} \longrightarrow Q$
Now, let $\nabla$ be a solution of the generalized HJ problem. Then one can find some solutions of the EL Equations, and the implicit Hamilton Equations, simply integrating $\nabla$.

We refer to Eq.{} (\ref{GHJE}) as \emph{Generalized HJ Equation} for the following reasons. Suppose that the Legendre transform is a diffeomorphism and search for those special solutions $\nabla$ of the generalized HJ problem such that
\begin{equation}
(\nabla , T)^\ast \Omega = 0, \label{dT}
\end{equation}
and
\begin{equation}
(\nabla, T)^\ast E. \label{HJE}
\end{equation}
Then, first of all, $\nabla$ is determined by $T$ via $\nabla := FL^{-1} \circ T$. Moreover, Eq.{} (\ref{dT}) is equivalent to $dT = 0$, so that, locally, $T = dS$ for some function $S$ on $Q$. Finally, Eq.{} (\ref{HJE}) is locally equivalent to $H \circ dS = 0$, $H$ being the Hamiltonian, which is the standard HJ Equation for $S$ (recall also that $\mathrm{im} T \subset T^\ast Q$ \emph{is a Lagrangian submanifold preserved by the dynamics}).
If $\nabla$ is locally given by $\nabla = \nabla^i \frac{\partial}{\partial q^i}$, then the generalized HJ Equation is locally given by
\begin{equation*}
\dfrac{\partial L}{\partial q^{i}}(q,\nabla)-\dfrac{\partial^{2}L}{\partial
\dot{q}^{i}\partial q^{j}}(q,\nabla)\nabla^{j}-\nabla^{k}\dfrac{\partial^{2}L}{\partial\dot
{q}^{i}\partial\dot{q}^{j}}(q,\nabla)\dfrac{\partial \nabla^{j}}{\partial q^{k}}  =0.
\end{equation*}

\section{Geometric Hamilton-Jacobi Field Theory}
We refer to \cite{b...99} for details on the geometry of jet spaces (see also \cite{v10,v11,v10b} for notations and conventions). Here, we only recall that the first jet bundle $\alpha_{1,0} : J^1 \alpha \longrightarrow N$ of a fiber bundle $\alpha : N \longrightarrow M$ is an affine bundle modelled over the vector bundle $\mathscr{V} := T^\ast M \otimes VN \longrightarrow N$, $VN \longrightarrow N$ being the vertical tangent bundle. The \emph{multimomentum bundle} of $\alpha$ is the bundle
\begin{equation*}
\alpha^\dag : J^\dag \alpha := \mathscr{V}^\ast \otimes \Lambda^n T^\ast M) \simeq V^\ast N \otimes \Lambda^{n-1} T^\ast M \longrightarrow M, \quad n = \dim M.
\end{equation*}
Local bundle coordinates $\ldots, x^i, \ldots, v^a, \ldots$ in $N$ determine standard coordinates $\ldots, x^i, \ldots, v^a,\ldots, v^a_i, \ldots$ in $J^1 \alpha$ and \textquotedblleft adjoint\textquotedblright{} coordinates $\ldots, x^i, \ldots, v^a,\ldots, p^i_a, \ldots$ in $J^\dag \alpha$.

Now, let $\pi : E \longrightarrow M $ be a configuration bundle, and $J^{k+1} \pi$ its $(k+1)$th jet bundle. Bundle coordinates $\ldots, x^i, \ldots, u^\alpha, \ldots$ in $E$ determine standard coordinates $\ldots,x^i,\ldots, u^\alpha_I, \ldots$, $|I|{} \leq k+1$ in $J^{k+1} \pi$. A \emph{Lagrangian field theory} of the order $k+1$ on $\pi$ is the datum of a Lagrangian density $\mathscr{L}$ on $J^{k+1} \pi$, i.e., a basic $n$-form on $J^{k+1} \pi$ (with respect to the projection $\pi_{k+1} : J^{k+1} \pi \longrightarrow M$). It determines a system of EL Equations of the order $2k + 2$, $\mathscr{E}_{EL} \subset J^{2k+2}$ for sections $\gamma : M \longrightarrow E$ of $\pi$. Roughly speaking, extremals of the action functional $\int \mathscr{L}$ are solutions of the EL equations, i.e., sections $\gamma$ of $\pi$ such that their jet prolongation $j_{2k+2}\gamma : M \longrightarrow J^{2k+2} \pi$ take values in $\mathscr{E}_{EL}$. Furthermore, $\mathscr{L}$ determines a system of first oder, implicit Hamilton-like PDEs $\mathscr{E}_H \subset J^1 \pi_k^\dag$ for sections of $\alpha^\dag$, as follows.

Put $\mathbb{J}^{k+1} \pi := J^{k+1} \pi \times_{J^k \pi} J^\dag \pi_k$, and let $p: \mathbb{J}^{k+1} \pi  \longrightarrow J^\dag \pi_k$ and $q : \mathbb{J}^{k+1} \longrightarrow J^{k+1} \pi$ be canonical projections. First of all, the Lagrangian density $\mathscr{L}$ determines a first order PDE for sections of $\mathbb{J}^{k+1} \longrightarrow M$. Namely, it holds the
\begin{proposition}
\cite{v10b} The Lagrangian density $\mathscr{L}$ determines an \emph{Hamiltonian} $(2,n-1)$-semibasic $(n+1)$-form (with respect to the projection $\mathbb{J}^{k+1} \longrightarrow M$) in a canonical way . If $\mathscr{L}$ is locally given by $\mathscr{L} = L d^n x$, $d^n x := dx^1 \wedge \cdots \wedge dx^n$, then $\Omega_\mathscr{L}$ is locally given by
\[
\Omega_{\mathscr{L}} = d p^{I.i}_\alpha \wedge d u^\alpha_I\wedge d^{n-1}x_i - dE \wedge d^n x, \quad
E = p_\alpha^{I.i} u^\alpha_{Ii} - L,
\]
with $d^{n-1} x_i := i (\partial / \partial x^i) d^n x$.
\end{proposition}
The first jet prolongation $j_1 \Sigma$ of a section $\Sigma$ of $\mathbb{J}^{k+1} \longrightarrow M$ can be \textquotedblleft inserted \textquotedblright{} into $\Omega_\mathscr{L}$ and the result $i(j_1 \Sigma) \Omega_\mathscr{L}$ is a section of a suitable vector bundle over $M$ (see \cite{v09b} for details). Therefore, we can consider the PDE
\begin{equation}
i({j_1 \Sigma}) \Omega_\mathscr{L} = 0. \label{ELHPDE}
\end{equation}
As usual, the geometric portrait of Eq.{} (\ref{ELHPDE}) is a subset $\mathscr{E}_{ELH}$ of the first jet bundle of $\mathbb{J}^{k+1} \longrightarrow M$. Since $\mathbb{J}^{k+1} $ projects onto $J^\dag \pi_k$, then $\mathscr{E}_{ELH}$ projects onto a suitable subset $\mathscr{E}_{H}$ of $J^1 \alpha^\dag$ which, in its turn, can be interpreted as a first order PDE for sections of $\alpha^\dag$. We refer to $\mathscr{E}_{H}$ as \emph{implicit Hamilton PDEs}.  

Eq.{} (\ref{ELHPDE}) covers both the EL Equations and the implicit Hamilton PDEs, i.e., if $\Sigma$ is a solution of Eq.{} (\ref{ELHPDE}), then $q \circ \Sigma = j_{k+1}{\gamma}$ for a solution $\gamma$ of the EL Equations, and $p \circ \Sigma $ is a solution of the implicit Hamilton PDEs. For this reason, I will refer to Eq.{} (\ref{ELHODE}) as \emph{ELH PDE}. Solutions of the ELH PDE take values in a suitable subbundle $\mathscr{P} \subset \mathbb{J}^{k+1}\pi$ of $q : \mathbb{J}^{k+1}\pi \longrightarrow J^{k+1} \pi$. Notice that the canonical projection $\mathscr{P} \longrightarrow J^k \pi$ is a bundle whose sections are pairs $(\nabla, T)$ such that $\nabla$ is a section of $\pi_{k+1,k}$ (in particular, an Ehresmann connection in $\pi_k$), and $T$ is a section of $J^\dag \pi \longrightarrow J^k \pi$.

Recall that the EL Equations have the local expression:
\begin{equation*}
(-)^{|I|}D_I \frac{\partial L}{\partial u^\alpha_ I} =0,
\end{equation*}
where the $D_I$'s are multiple total derivatives. The ELH PDEs have the following local expression:
\[
\left\{ 
\begin{array}{ll}
\dfrac{\partial }{\partial x^{i}}u_{I}^{\alpha }=u_{Ii}^{\alpha } & \quad
|I|{}\leq k \\ 
\vspace{-0.2cm} &  \\ 
\dfrac{\partial }{\partial x^{i}}p_{\alpha }^{I.i}=\dfrac{\partial L}{%
\partial u_{I}^{\alpha }}[x,u]-\delta _{Ji}^{I}p_{\alpha }^{J.i} & \quad |I|{}\leq
k \\ 
\vspace{-0.2cm} &  \\ 
\delta _{Ji}^{I}p_{\alpha }^{J.i}=\dfrac{\partial L}{\partial u_{I}^{\alpha }%
}[x,u] & \quad |I|{}=k+1%
\end{array}%
\right. ,
\]
where, $\delta_K^I$ is equal to $1$ if $I = K$ and is equal to $0$ otherwise.
In particular $\mathscr{P}$ is locally given by 
\[
\delta _{Ji}^{I}p_{\alpha }^{J.i}=\dfrac{\partial L}{\partial u_{I}^{\alpha }%
}[x,u], \quad |I|{} = k+1.
\]
The implicit Hamilton PDEs have the local expression
\[
\left\{ 
\begin{array}{ll}
\dfrac{\partial }{\partial x^{i}}u_{I}^{\alpha }=u_{Ii}^{\alpha } & \quad |I|{}<k
\\ 
\vspace{-0.2cm} &  \\ 
\dfrac{\partial }{\partial x^{i}}p_{\alpha }^{I.i}=\dfrac{\partial L}{%
\partial u_{I}^{\alpha }}[x,\bar{u}]-\delta _{Ji}^{I}p_{\alpha }^{J.i} & \quad
|I|{}\leq k \\ 
\vspace{-0.2cm} &  \\ 
\delta _{Ji}^{I}p_{\alpha }^{J.i}=\dfrac{\partial L}{\partial u_{I}^{\alpha }%
}[x,\bar{u}] & \quad |I|{}=k+1%
\end{array}%
\right. ,
\]%
where $\bar{u}=(\ldots ,\bar{u}_{I}^{\alpha },\ldots )$, $|I|{}\leq k+1$, $\bar{u}_{I}^{\alpha }: = {u}_{I}^{\alpha }$ for $|I|{} \leq k$, while $ \bar{u}_{Ii}^{\alpha } := \partial u^\alpha_I / \partial x^i$ for $|Ii| = k+1$.

The central problems in the geometric HJ field theory hare the following

\begin{problem}[Lagrangian Generalized HJ Problem]
Find a flat, holonomic connection $\nabla$ in $\pi_ k$ (see \cite{v11}) such that any $\nabla$-constant section $\Gamma : M \longrightarrow J^k \pi$ is of the form $\Gamma = j_k \gamma$, where $\gamma$ is a solution of the EL Equations.
\end{problem}

\begin{problem}[Hamiltonian Generalized HJ Problem]
Find a flat, holonomic connection $\nabla$ in $\pi_ k$ and a section $T$ of $J^\dag \pi_k \longrightarrow J^k \pi$ such that, for any $\nabla$-constant section $\Gamma : M \longrightarrow J^k \pi$, $T \circ \Gamma$ is a solution of the implicit Hamilton PDEs.
\end{problem}

\begin{problem}[Generalized HJ Problem]
Find a flat, holonomic connection $\nabla$ in $\pi_ k$ and a section $T$ of $J^\dag \pi_k \longrightarrow J^k \pi$ such that, for any $\nabla$-constant section $\Gamma : M \longrightarrow J^k \pi$, $(\nabla, T) \circ \Gamma$ is a solution of the ELH PDEs.
\end{problem}

Before stating the key result in the geometric HJ field theory recall that a connection $\nabla$ in $\pi_k$ can be \textquotedblleft inserted\textquotedblright{} into a $(2, \ell)$-semibasic $(2+\ell)$-form $\omega$ on $J^k \pi$ and the result $i(\nabla) \omega$ is a section of a suitable vector bundle over $J^k \pi$ (see \cite{v09b} for details).

\begin{theorem}[Generalized HJ Theorem]
Let $\nabla$ be a flat, holonomic connection in $\pi_ k$. The following conditions are equivalent:
\begin{itemize}
\item $\nabla$ is a solution of the Lagrangian generalized HJ problem;
\item $(\nabla, T)$ is a solution of the Hamiltonian generalized HJ problem for some section $T$ of $J^\dag \pi_k \longrightarrow J^k \pi$.
\item $(\nabla, T)$ is a solution of the generalized HJ problem for some section $T$ of $J^\dag \pi_k \longrightarrow J^k \pi$.
\item $i(\nabla) (\nabla , T)^\ast \Omega_\mathscr{L} = 0$, for some section $(\nabla , T)$ of $\mathscr{P} \longrightarrow J^k \pi$.
\end{itemize}
\end{theorem}
In particular, solving the generalized HJ problem amounts to solving the PDE 
\begin{equation}
i(\nabla) (\nabla , T)^\ast \Omega_\mathscr{L} = 0, \label{GFHJE}
\end{equation}
Now, let $(\nabla, T)$ be a solution of the generalized HJ problem. Then one can find some solutions of the EL Equations and the implicit Hamilton PDEs, simply integrating $\nabla$ (Examples can be found in \cite{v10,v11}).

Since Eq.{} (\ref{GFHJE}) generalizes Eq.{} (\ref{GHJE}) to the field theoretic setting, we again refer to it as \emph{Generalized HJ Equation}. The main difference between the mechanical and the field theoretic cases is that in the latter, $T$ is not uniquely determined by $\nabla$ and, therefore, the Lagrangian and the Hamiltonian problems are not strictly equivalent. This is a consequence of the fact that there is no uniquely defined Legendre transform in higher derivative field theory. Nonetheless, $T$ is uniquely determined by $\nabla$ (and $\mathscr{L}$) up to \textquotedblleft total divergences\textquotedblright. As a corollary the generalized HJ problem does only depend on the EL Equations (and not on the Lagrangian density) up to isomorphisms (see \cite{v10,v11} for details).
Finally notice that a natural question would be the following: how does the standard HJ problem generalize to the field theoretic setting? A natural candidate for the field theoretic version of the standard HJ Equation is $(\nabla, T) ^\ast \Omega_\mathscr{L} = 0$. However, this answer is unsatisfactory (see \cite{v11}). In fact, there are indications that a better field theoretic version of the standard HJ theory, encompassing all its aspects (including the HJ Equations, its complete integral, canonical transformations, etc.) lives on the \emph{diffiety of initial data of the EL Equations}. Work on this conjecture is still in progress.

We conclude this report summarizing the theory in a table of correspondences between mechanical concepts and their field theoretic analogues.

\begin{center}
\begin{tabular}
[c]{cc|cc}%
\bf{Mechanics} &  &  & \bf{Field Theory} \\
\hline\hline
{$Q$}${}\times \{ \text{time} \}$ &  &  & {$J^k \pi$}\\

{$TQ$}${}\times \{ \text{time} \}$ &  &  & {$J^{k+1} \pi$}\\

{$T^2 Q$}${}\times \{ \text{time} \}$ &  &  & {$J^{2k+2} \pi$}\\

{$T^{\ast}Q$}${}\times \{ \text{time} \}$ &  &  & {$J^{\dag}\pi_k$} \\

{$\mathbb{T}Q$}${}\times \{ \text{time} \}$ &  &  & {$\mathbb{J}^{k+1} \pi$} \\

{$\Omega - dE \wedge dt$} & & & {$\Omega_\mathscr{L} $} \\

{graph of $FL$} & & & {$\mathscr{P}$} \\

{$i({\dot{\Sigma}}) \Omega = dE |_{\Sigma}$} & & & $i({j_1 \Sigma}) \Omega = 0$ \\

$i(\nabla) (\nabla, T)^\ast \Omega = (\nabla, T)^\ast E$ & & & $i(\nabla) (\nabla, T)^\ast \Omega_\mathscr{L} = 0$

\end{tabular}

\end{center}

\section*{Acknowledgements}

I thank the organizers of the workshop \emph{FunInGeo, 08--12 June 2011, Ischia (NA) Italy, in honour of Giuseppe Marmo's 65th birthday}, for the opportunity of presenting this material at the meeting. This short report is dedicated to Beppe: \emph{he is the one who directed me to Geometry}.

\end{document}